\documentclass{article}
\usepackage{ast2}
\title{W-types in sheaves}
\author{Benno van den Berg \& Ieke Moerdijk}
\date{September 25, 2008}
\begin{document}

\maketitle

\begin{abstract}
\noindent
In this small note we give a concrete description of W-types in categories of sheaves.
\end{abstract}

\noindent
It can be shown that any topos with a natural numbers object has all W-types. Although there is this general result, it can be useful to have a concrete description of W-types in various toposes. For example, a concrete description of W-types in the effective topos can be found in \cite{berg06,hofmannvanoostenstreicher06}, and a concrete description of W-types in categories of presheaves was given in \cite{moerdijkpalmgren00}. It was claimed in \cite{moerdijkpalmgren00} that W-types in categories of sheaves are computed as in presheaves (Proposition 5.7 in \emph{loc.cit.}) and can therefore be described in the same way. Unfortunately, this claim is incorrect, as the following (easy) counterexample shows. Let $f: 1 \to 1$ be the identity map on the terminal object. The W-type associated to $f$ is the initial object, which, in general, is different in categories of presheaves and sheaves. This means that we still lack a concrete description of W-types in categories of sheaves. This note aims to fill this gap.

We would like to warn readers who are sensitive to such issues that our metatheory is {\bf ZFC}. \emph{In particular, we freely use the axiom of choice}. We leave the issue of how to describe W-types in categories of sheaves when the metatheory is more demanding (i.e. weaker) to another occasion.

Categories of sheaves are described using (Grothendieck) sites. There are different formulations of the notion of a site, all essentially equivalent (\cite{johnstone02b} provides an excellent discussion of this point), but for our purposes we find the following (``sifted'') formulation the most useful. 

\begin{defi}{site}
Let \ct{C} be a category. A \emph{sieve} $S$ on an object $a \in \ct{C}$ consists of a set of arrows in \ct{C} all having codomain $a$ and closed under precomposition (i.e., if $f: b \to a$ and $g: c \to b$ are arrows in \ct{C} and $f$ belongs to $S$, then so does $fg$). We call the set $M_a$ of all arrows into $a$ the \emph{maximal sieve} on $a$. If $S$ is a sieve on $a$ and $f: b \to a$ is any map in \ct{C}, we write $f^*S$ for the sieve $\{ g: c \to b \, : \, fg \in S \}$ on $b$. In case $f$ belongs to $S$, we have $f^*S = M_b$.

A \emph{(Grothendieck) topology} Cov on \ct{C} is given by assigning to every object $a \in \ct{C}$ a collection of sieves ${\rm Cov}(a)$ such that the following axioms are satisfied:
\begin{description}
\item[(Maximal sieve)] The maximal sieve $M_a$ belongs to ${\rm Cov}(a)$;
\item[(Stability)] If $f: b \to a$ is any map and $S$ belongs to ${\rm Cov}(a)$, then $f^*S$ belongs to ${\rm Cov}(b)$;
\item[(Local character)] If $S$ is a sieve on $a$ for which there can be found a sieve $R \in {\rm Cov}(a)$ such that for all $f: b \to a \in R$ the sieve $f^*S$ belongs to ${\rm Cov}(b)$, then $S$ belongs to ${\rm Cov}(a)$. 
\end{description}
A pair $(\ct{C}, {\rm Cov})$ consisting of a category $\ct{C}$ and a topology Cov on it is called a \emph{site}. If a site $(\ct{C}, {\rm Cov})$ has been fixed, we call the sieves belonging to some ${\rm Cov}(a)$ \emph{covering sieves}. If $S$ belongs to ${\rm Cov}(a)$ we say that $S$ is a \emph{sieve covering} $a$, or that $a$ \emph{is covered by} $S$. 
\end{defi}

With our notion of site in place, we are ready to describe W-types in sheaves. First, we fix a map $F: Y \to X$ of sheaves. Let ${\cal V}$ be the set of all well-founded trees, in which 
\begin{itemize}
\item nodes are labelled with triples $(a, x, S)$ with $a$ an object in \ct{C}, $x$ an element of $X(a)$ and $S$ a sieve covering $a$,
\item edges are labelled with pairs $(f, y)$ with $f: b \to a$ a map in \ct{C} and $y$ an element of $Y(b)$,
\end{itemize}
in such a way that
\begin{itemize}
\item if a node is labelled with $(a, x, S)$ and an edge into this node is labelled with $(f, y)$, then $f$ belongs to $S$ and $F(y) = x \cdot f$,
\item and if a node is labelled with $(a, x, S)$ there is, for every $f \in S$ and $y$ such that $F(y) = x \cdot f$, a unique edge into this node labelled with $(f, y)$.
\end{itemize}
In fact, ${\cal V}$ is the W-type (in the category of sets) associated to the projection
\diag{  \sum_{(a, x, S)} \{ (f, y) \, : \, f: b \to a \in S, y \in F(b), F(y) = x \cdot f \} \ar[d] \\ \{ (a, x, S) \, : \, a \in \ct{C}, x \in X(a), S \in {\rm Cov}(a) \}.}
If $v$ denotes a well-founded tree in ${\cal V}$, we will also use the letter $v$ for the function that assigns to labels of edges into the root of $v$ the tree attached to this edge. So if $(f, y)$ is a label of one of the edges into the root of $v$, we will write $v(f, y)$ for the tree that is attached to this edge; this is again an element of ${\cal V}$. Note that an element of  ${\cal V}$ is uniquely determined by the label of its root and the function we just described. 

We say that a tree $v \in {\cal V}$ is \emph{rooted} at an object $a$ in \ct{C}, if its root has a label whose first component is $a$. We will denote the collection of trees rooted at $a$ by ${\cal V}(a)$. This gives the set ${\cal V}$ the structure of a presheaf, as one can define the following restriction operation. Let $v \in {\cal V}(a)$ and $f: b \to a$ be a map in \ct{C}. If $v$ has root $(a, x, S)$ the tree $v \cdot f$ has root $(b, x \cdot f, f^* S)$ and
\[ (v \cdot f)(g, y) = v(fg, y). \]
One easily verifies that this is well-defined, and gives $V$ the structure of a presheaf.

Next, we define by transfinite recursion an equivalence relation $\sim$ on the presheaf ${\cal V}$:
\begin{center}
\begin{tabular}{lcp{8 cm}}
$v \sim v'$ & $\Leftrightarrow$ & If the root of $v$ is labelled with $(a, x, S)$ and the root of $v'$ with $(a', x', S')$, then $a = a'$, $x = x'$ and there is a covering sieve $R \subseteq S \cap S'$ such that for every $f: b \to a \in R$ and $y \in Y(b)$ such that $F(y) = x \cdot f$ we have $v(f, y) \sim v'(f, y)$. 
\end{tabular}
\end{center}
By transfinite induction one verifies that $\sim$ is an equivalence relation. Furthermore, one verifies directly that $\sim$ is a presheaf (i.e. $v \sim v'$ implies that $v$ and $v'$ are rooted at the same object $a$ and that $v \cdot f \sim v' \cdot f$ for all $f: b \to a$).

Next, we define \emph{composability} and \emph{naturality} of trees (the terminology is taken from \cite{moerdijkpalmgren00}):
\begin{itemize}
\item A tree $v$ whose root is labelled with $(a, x, S)$ is \emph{composable}, if for any $(f, y)$ with $f: b \to a \in S$ and $y \in Y(b)$ such that $F(y) = x \cdot f$, the tree $v(f, y)$ is rooted at $b$.
\item A tree $v$ whose root is labelled with $(a, x, S)$ is \emph{natural}, if for any $(f, y)$ with $f: b \to a \in S$, $g: c \to b$ and $y \in Y(b)$ such that $F(y) = x \cdot f$, we have
\[ v(f, y) \cdot g \sim v(fg, y \cdot g). \]
\end{itemize}
It is clear that natural and composable trees are stable under restriction, so that also these form presheaves. The same applies to the presheaf ${\cal W}$ of trees that are \emph{hereditarily} composable and natural (i.e. not only are they themselves both composable and natural, but the same is true for all their subtrees). 

The relation $\sim$ is also an equivalence relation on ${\cal W}$ and we let $\overline{\cal W}$ be its quotient \emph{in presheaves} (so the quotient is computed pointwise). We show that $\overline{\cal W}$ is a sheaf and, indeed, the W-type associated to $F$ in sheaves.

\begin{lemm}{Wbarseparated}
If $T$ is a sieve covering $a$ and $w, w' \in {\cal W}(a)$ are such that $w \cdot f \sim w' \cdot f$ for all $f \in T$, then $w \sim w'$. In other words, $\overline{\cal W}$ is separated.
\end{lemm}
\begin{proof} If the label of the root of $w$ is of the form $(a, x, S)$ and that of $w'$ is of the form $(a, x', S')$, then $w \cdot f \sim w' \cdot f$ implies that $x \cdot f = x' \cdot f$ for all $f \in T$. As $X$ is separated, it follows that $x = x'$.

Consider
\begin{eqnarray*}
R & = & \{ \, g: b \to a  \in (S \cap S') \, : \, \forall h: c \to b, y \in Y(c) \, \\
& & [ \, F(y) = x \cdot gh \Rightarrow w(gh, y) \sim w'(gh, y)\, ] \, \}.
\end{eqnarray*}
$R$ is a sieve, and the statement of the lemma will follow once we show that it is covering.

Fix an element $f \in T$. That $w \cdot f \sim w' \cdot f$ holds means that there is a covering sieve $R_f \subseteq f^*S \cap f^*S'$ such that for every $k: c \to b \in R_f$ and $y \in Y(c)$ such that $F(y) = x \cdot fk$ we have $w(fk, y) = (w \cdot f)(k, y) \sim (w' \cdot k)(g, y) = w'(fk, y)$. In other words, $R_f \subseteq f^*R$. So $R$ is a covering sieve by local character.
\end{proof}

\begin{lemm}{Wbarsheaf}
$\overline{\cal W}$ is a sheaf.
\end{lemm}
\begin{proof} Let $S$ be a covering sieve on $a$ and suppose we have a compatible family of elements $(\overline{w}_f \in \overline{\cal W})_{f \in S}$. For every element $f \in S$ \emph{choose} a representative $(w_f \in {\cal W})_{f \in S}$ such that $[w_f] = \overline{w}_f$. For every $f: b \to a \in S$ a representative $w_f$ has a root labelled by something of form $(b, x_f, R_f)$. The $x_f$ form a compatible family and, since $X$ is a sheaf, can be glued together to obtain an element $x \in X(a)$. Furthermore, $R := \{ fg \, : \, f \in S, g \in R_f \} \subseteq S$ is a covering sieve, by local character.

Before we proceed, we prove the following claim: \\
\noindent
{\bf Claim.} Assume $f: b \to a \in S$ and $g: c \to b \in R_f$ and $f': b' \to a \in S$ and $g': c \to b' \in R_{f'}$ are such that $fg = f'g'$. If $y \in Y(c)$ is such that $F(y) = x \cdot fg$, then
\[ w_f(g, y) \sim w_{f'}(g', y). \]
\begin{proof} By compatibility of the family $(\overline{w}_f \in \overline{\cal W})_{f \in S}$ we know that $w_f \cdot g \sim w_{f'} \cdot g' \in {\cal W}(c)$. This means that there is a covering sieve $T \subseteq g^*R_f \cap (g')^*R_{f'}$ such that for all $h: d \to c \in T$ and $z \in Y(d)$ such that $F(z) = x \cdot fgh$, we have $(w_f \cdot g)(h, z) \sim (w_{f'} \cdot g')(h, z)$. So if $h: d \to c \in T$, then  
\begin{eqnarray*}
w_f(g, y) \cdot h & \sim & w_f(gh, y \cdot h) \\
& = & (w_f \cdot g)(h, y \cdot h) \\
& \sim & (w_{f'} \cdot g')(h, y \cdot h) \\
& = & w_{f'}(g'h, y \cdot h) \\
& \sim & w_{f'}(g', y) \cdot h.
\end{eqnarray*}
Because $\overline{\cal W}$ is separated (as was shown in \reflemm{Wbarseparated}), it follows that $w_f(g, y) \sim w_{f'}(g', y)$. This completes the proof of the claim.
\end{proof}

We construct a tree $w \in {\cal V}$ such that $[w]$ is the desired glueing of all the $\overline{w}_f$. It will have a root labelled with $(a, x, R)$. For any $h: c \to a \in R$ and $y \in Y(c)$ such that $F(y) = x \cdot h$, we \emph{choose} $f \in S$ and $g \in R_f$ such that $h = fg$ and set
\[ w(h, y) = w_f(g, y). \]
This \emph{does} essentially depend on the choice of $f$ and $g$, but any two choices yield \emph{equivalent} results: that is precisely the content of the claim we proved above.

It is easy to see that the tree that we have constructed is composable. It is also natural, since if $h: c \to a \in R$ and $y \in Y(c)$ are such that $F(y) = x \cdot h$, and we have chosen $f \in S$ and $g \in R_f$ such that $h = fg$ and we have set
\[ w(h, y) = w_f(g, y), \]
and $k: d \to c$ is any other map, and we have chosen $f' \in S$ and $g' \in R_{f'}$ such that $hk = f'g'$ and we have set
\[ w(hk, y \cdot k) = w_{f'}(g', y \cdot k), \]
then it follows that
\begin{eqnarray*}
w(h, y) \cdot k & = & w_f(g,y) \cdot k \\
& \sim & w_f(gk, y \cdot k) \\
& \sim & w_{f'}(g', y \cdot k) \qquad \mbox{(using the claim)} \\
& = & w(hk, y \cdot k),
\end{eqnarray*}
as desired.

It remains to show that $[w]$ is a glueing of all the $\overline{w}_f$, i.e. that $w \cdot f \sim w_f$ for all $f \in S$. So let $f: b \to a \in S$. First of all, $x \cdot f = x_f$, by construction. Furthermore, for every $g: c \to b \in R_f = R_f \cap f^*R$ and $y \in Y(c)$ such that $F(y) = x \cdot fg$, let $f' \in S$ and $g' \in R_{f'}$ be such that $fg = f'g'$ and $w(fg, y) = w_{f'}(g', y)$. Then
\begin{eqnarray*}
(w \cdot f)(g, y) & = & w(fg, y) \\
& = & w_{f'}(g', y) \\
& \sim & w_{f}(g, y) \qquad \mbox{(using the claim)}.
\end{eqnarray*}
This completes the proof.
\end{proof}

\begin{lemm}{Wbaralgebra}
$\overline{\cal W}$ is a $P_F$-algebra.
\end{lemm}
\begin{proof} We have to describe a natural transformation ${\rm sup}: P_F \overline{\cal W} \to \overline{\cal W}$. An element of $P_F \overline{\cal W}(a)$ is a pair $(x, t)$ consisting of an element $x \in X(a)$ together with a natural transformation $t: Y_x \to \overline{\cal W}$, where $Y_x$ is the presheaf defined by
\[ Y_x(b) = \{ (f: b \to a, y) \, : \, F(y) = x \cdot f \}. \]
We define ${\rm sup}_x(t)$ to be $[w]$, where tree $w$ is the tree whose root has label $(a, x, M_a)$ and for which, for every $(f, y) \in Y_x$, the value of $w(f, y)$ is \emph{chosen} such that $[w(f, y)] = t(f, y)$ (this is another application of choice). One quickly verifies that $w$ is composable and natural.

The other verification (that of the naturality of the {\rm sup}-operation) is easy and also left to the reader.
\end{proof}

\begin{lemm}{Wbaralgebra}
$\overline{\cal W}$ is the initial $P_F$-algebra.
\end{lemm}
\begin{proof}
We follow the usual strategy: we show that ${\rm sup}: P_F\overline{\cal W} \to \overline{\cal W}$ is monic and that $\overline{\cal W}$ has no proper $P_F$-subalgebras (i.e., we apply Theorem 26 of \cite{berg05}).

We first show that sup is monic. So let $(x, t), (x', t') \in P_FX(a)$ be such that ${\rm sup}_x(t)={\rm sup}_{x'}(t')$. It follows that $x = x'$ and that there is a covering sieve $S$ on $a$ such that for all $h: b \to a \in S$ and $y \in Y(b)$ such that $F(y) = x \cdot h$, we have $t(h, y) = t'(h,  y)$. We need to show that $t = t'$, so let $(f, y) \in Y_x$ be arbitrary. For every $g \in f^*S$, we have:
\begin{eqnarray*}
t(f, y) \cdot g & = & t(fg, y \cdot g) \\
& = & t'(fg, y \cdot g) \\
& = & t'(f, y) \cdot g.
\end{eqnarray*}
Since $f^*S$ is covering, it follows that $t(f, y) = t'(f, y)$, as desired.

The fact that $\overline{\cal W}$ has no proper $P_F$-subalgebras is a consequence of the inductive properties of ${\cal V}$ (we recall that ${\cal V}$ is a W-type). Let $I$ be a sheaf and $P_F$-subalgebra of $\overline{\cal W}$. We claim that 
\[ J = \{ v \in {\cal V} \, : \, \mbox{if } v \mbox{ is both hereditarily composable and natural, then } [v] \in I \} \]
is a subalgebra of ${\cal V}$. Proof: Suppose $v$ is a tree that is both hereditarily composable and natural. Assume moreover that $(a, x, S)$ is the label of its root and that for all $(f, y)$ with $f: b \to a \in S$ and $F(y) = x \cdot f$, we know that $[v(f, y)] \in I$. Our aim is to show that $[v] \in I$.

For the moment fix an element $f: b \to a \in S$. The tree $v \cdot f$ has $(b, x \cdot f, M_b)$ as the label of its root and for any $(g, y)$ with $g: c \to b \in M_b$ and $F(y) = x \cdot fg$ the tree $(v \cdot f)(g, y)$ is given by $v(fg, y)$. This means that $[ v] \cdot f  = {\rm sup}_{x \cdot f}(t)$, where $t(g, y) = [v(fg, y)] \in I$. Since $I$ is a $P_F$-subalgebra of $\overline{\cal W}$ this implies that $[v] \cdot f \in I$. From this it follows that $[v] \in I$, since $I$ is also a sheaf and $f$ was an arbitrary element of the covering sieve $S$.

We conclude that $J = {\cal V}$ and $I = \overline{\cal W}$. This completes the proof.
\end{proof}

This completes the proof of the correctness of our description of W-types in categories of sheaves.  

\newpage

\bibliographystyle{plain} \bibliography{ast}

\end{document}